\documentclass[twoside,leqno,fleqn]{article}
\usepackage{amsmath}
\usepackage{amssymb}

\usepackage{calligra}
\usepackage[T1]{fontenc}

\usepackage[shortlabels]{enumitem}
\usepackage[onehalfspacing]{setspace}

\newcommand{\AUTOR}{C. Ga\ss ner}
\newcommand{\TITEL}{AC and the Independence of WO in Second-Order Henkin Logic, Part II}
\def\S{\Sigma}
\def\I{{\cal I}}

\def\G{\mathfrak{G}}
\def\F{\mbox{\scriptsize\calligra T\!\!\!\!\!\hspace{0.005cm}F\!\!\!\!\!T\!\!\!\!\!\hspace{0.02cm}F}\,\,\,}

\newcommand{\bbbn}{\mathbb{N}} 

\newcommand{\mbm}[1]{\mbox{\boldmath{$#1$}}}
\newcommand{\mbmss}[1]{\mbox{\scriptsize\boldmath{$#1$}}}
\newcommand{\mbmty}[1]{\mbox{\tiny\boldmath{$#1$}}}
\newcommand{\qed}{\hfill{$\Box$}} 

\newtheorem{satz}{Satz}[section] 
\newtheorem{lemma}[satz]{Lemma}
\newtheorem{proposition}[satz]{Proposition}
\newtheorem{corollary}[satz]{Corollary}
\newtheorem{theorem}[satz]{Theorem}

\newtheorem{agree}[satz]{Agreement}

\begin{document}

\newcounter{li}

\thispagestyle{empty}
\begin{center} {\Large\bf AC and the Independence of WO \vspace{0.4cm}\\ in Second-Order Henkin Logic, Part II\,\footnote{\footnotesize For Part I see arXiv:2409.10276 [math.LO].}}\vspace{0.6cm}\\{\bf Christine Ga\ss ner}\footnote{I thank Michael Rathjen for the discussion on various questions related to the paper of Siskind, Mancosu, and Shapiro (2020) in Leeds. My thanks also go to the organizers of the Colloquium Logicum in Konstanz in 2022. } {\vspace{0.2cm}\\University of Greifswald, Germany, 2024\\ gassnerc@uni-greifswald.de}\\\end{center} 

\begin{abstract} This article is concerned with the Axiom of Choice (AC) and the well-ordering theorem (WO) in second-order predicate logic with Henkin interpretation (HPL). We consider a principle of choice introduced by Wilhelm Ackermann (1935) and discussed also by David Hilbert and Ackermann (1938), by Günter Asser (1981), and by Benjamin Siskind, Paolo Mancosu, \linebreak and Stewart Shapiro (2020). 
Our discussion is restricted to so-called Henkin-Asser structures of second order.  Here, we give the technical details of our proof of the independence of WO from the so-called Ackermann axioms in HPL presented at the Colloquium Logicum in 2022. Most of the definitions used here can be found in Sections 1, 2, and 3 of Part I. 

\end{abstract}

\setcounter{section}{3}

\section{The independence of the second-order WO}\label{SectionIndependWO}

In this section, we consider the basic Fraenkel model $\S_0$ of second order which we also call the  basic Fraenkel model II for short.  It is defined by  $\S_0= \S(\bbbn, \G^\bbbn_{\sf 1},\F_0(\bbbn, \G^\bbbn_{\sf 1})) $ where $\bbbn$ can be an arbitrary  infinite set. To simplify matters, we assume that ZFC is consistent and we use it as metatheory. Consequently,  $\bbbn$ is a well-ordered set  in our metatheory as usual.

Any structure $ \S(I, \G,\F_0(I, \G)) $ is determined by a subgroup $\G$ of $\G^I_{\sf 1}$ and the normal ideal $\I_0^I$ where $ \G^I_{\sf 1}$ is the group of all permutations in $perm(I)$ and $\I_0^I$ is the system of all finite subsets of $I$. We denote this structure also by $ \S(I, \G,\I_0^I)$ and by $(J_n(I, \G,\I_0^I))_{n\geq 0}$.\footnote{For details see Sections 1, 2, and 3 in Part I where we use the references that are listed again at the end.} 

\subsection{The well-ordering principle and the Ackermann axiom}\label{SectionWOandMore}

The well-ordering theorem $WO^n$ states that, for any $n$-ary predicate, there is a $2n$-ary well-ordering. In this section, we want to show that $WO^n$ cannot be deduced from the set $choice_h^{n,m}$ of Ackermann axioms in HPL. Let the formula $wo(T,A)$ describe the properties of a well-ordering $\widetilde\tau$ on a set $\widetilde{\alpha}$. Here, $T$ stands for the $2n$-ary predicate $\tau$ and $A$ for the $n$-ary predicate $\alpha$. Such a predicate $\alpha$ is {\em well-ordered} by the predicate $\tau$ if $\widetilde{\alpha}$ is well-ordered by $\widetilde{\tau}$. Therefore, $\alpha$ is well-ordered by $\tau$ if the set $\widetilde{\alpha}$ is linearly ordered by $\widetilde\tau$ and if each of its non-empty subsets $\widetilde{\beta}\subseteq \widetilde{\alpha}$ has a minimal element $\mbm{\xi}_0$ with respect to $\tau$.
More precisely, let $po(T,A)$ stand for the statement that $T$ is a partial order in $A$ and let

$lo(T,A)=_{\rm df} po(T,A)\land \,\forall \mbm{x}_1\forall \mbm{x}_2
 ( A\mbm{x}_1 \land A\mbm{x}_2 \to T\mbm{x}_1\mbm{x}_2
 \lor T\mbm{x}_2\mbm{x}_1)$. \label{DefLinOrd}

\vspace{0.2cm}

\noindent By analogy to \cite{1c}, we say that each linear ordering $\tau$ on a predicate $\alpha$ and each well-ordering are reflexive. Let

\vspace{0.2cm}
$wo(T,A)=_{\rm df}lo(T,A)\land\forall B ( \forall \mbm{x}(B\mbm{x}\!\to\! A\mbm{x})\land \exists\mbm{x}B\mbm{x}$

\hspace{5cm}$ \to \exists\mbm{x}_0(B\mbm{x}_0 \land \forall\mbm{x}_1(B\mbm{x}_1 \to T\mbm{x}_0\mbm{x}_1)))$. \label{DefWellOrd}

\vspace{0.2cm}
\noindent Consequently, let $WO^n$ be defined by 
\[WO^{n} =_{\rm df}\forall A \exists T wo(T,A).\]

By \cite{1c}, we have the following. 
\begin{proposition}[Asser {\normalfont \cite{1c}}]\label{PropertiesWO} For $n,m \geq 1$, there hold the following relationships. 
\[\models_hWO^1\leftrightarrow WO^n\]\[
\models_h WO^1\to AC^{n, m}\]
\end{proposition}
In 1981, the question whether all formulas in $choice_h$ are true in all models of $^{h}ax^{(2)}\cup\{WO^1\}$ seemed to be an unsolved problem and Asser \cite{1c} wrote that no known proof of $ WO^n $ can be transferred from ZFC to HPL if the axiom of choice is given in the form of formulas resulting from $choice_h$. Since 2018, these problems were discussed for $choice_h^{1,1}$
by Siskind, Mancosu, and Shapiro (cf. \cite{Siskind}).

We also will deal with Asser's second problem and show that $ WO^n $ cannot be deduced from $^{h}ax^{(2)}\cup choice_h$ (cf.\,\,Proposition \ref{WO_unabh_choice}).

\begin{proposition}\label{WOnichtInS1}For each $n \geq 1$, there holds 
\[\begin{array}{l}\S_0\models \neg WO^n.
\end{array}\]
\end{proposition}
{\bf Proof.} Let $LO^{n} =_{\rm df}\forall A \exists T lo(T,A)$. Then, the assumption follows from Proposition \ref{PropertiesWO} 
and the following facts, (a) and (b).

\vspace{0.2cm}
\noindent \begin{tabular}{lll}
 (a) & $\models_h WO^1 \to LO^1$ & \cite[p.\,28, Section IV, Satz 2.1.1, c)]{Gass84}\\ \\
(b) &$\S_0\models \neg LO^1$ & \cite[p.\,49, Section V, Satz 1.5]{Gass84} (see also \cite{Gass94})\hspace{0.6cm}\qed\\\end{tabular}

\vspace{0.4cm}
\noindent Consequently, to prove that $WO^1$ does not result from $choice_h$ in {\rm HPL}, it is sufficient (for any $n,m\geq 1$) to show that $choice_h^{n,m}(H)$ holds for all formulas $H\in {\cal L}^{(2)}_{\mbm x,D}$ in $\S_0$. 

Roughly speaking, the main idea for proving $\S_0\models choice_h^{n,m}(H)$ for each formula $H\in {\cal L}^{(2)}_{\mbm x,D}$ is to show the existence of the required predicate by decomposing the problem into subproblems so that it will be possible and sufficient to construct suitable sub-predicates in various ways. 
For this reason, we start with an analysis of the antecedent $\forall \mbm{x} \exists D H(\mbm{x},D)$ and the consequent $\exists S \forall \mbm{x} \exists D (\forall\mbm{y}(D \mbm{y} \leftrightarrow S\mbm{x}\mbm{y}) \land H(\mbm{x}, D ))$ of the Ackermann axiom $choice^{n,m}_h(H)$ for any second-order formula $H(\mbm{x},D)$ in ${\cal L} ^{(2)}_{\mbmss{x},D}$. 

\begin{lemma}[The assumption $\forall \mbm{x} \exists D H(\mbm{x},D)$ and choice functions]\label{Zhg_CG_sigma}\hfill Let\linebreak $\S$ be any predicate structure $(J_n)_{n\geq 0}$ in ${\sf struc}_{\rm pred}^{({\rm m})}$, $f$ be any assignment in $\S$, and $n,m\geq 1$. For any  formula $H$  in ${\cal L} ^{(2)}_{\mbmss{x},D}$, let 
\[ {\cal A}^H=\{(\mbm{\xi}, \delta )\in J_0^n\times J_m\mid \S\models_{f\langle{\mbmty{x}\atop \mbmty{\xi}}{D\atop\delta}\rangle} H\}\]
 and, for any $\mbm{\xi}\in J_0^n$, let
\[{\cal A}^H_{\mbmss{\xi}}=\{ \delta \in J_m\mid (\mbm{\xi}, \delta )\in {\cal A}^H\}.\] 

\begin{enumerate}[label={\em(\arabic*)}]
 \labelwidth0.7cm \leftmargin0cm \itemsep2pt plus1pt
\topsep1pt plus1pt minus1pt
\labelsep4pt \parsep0.5pt plus0.1pt minus0.1pt
\item $\S\models_f\forall \mbm{x} \exists D ( H(\mbm{x},D))$ implies that, for all lists $\mbm{\xi}$ of $n$ individuals in $J_0$, there is a non-empty uniquely determined subset ${\cal A}^H_{\mbmss{\xi}}\subseteq J_m$.
\item If there is a choice function $\varphi$ on the family $({\cal A}^H_{\mbmss{\xi}})_{\mbmss{\xi}\in J_0^n}$ such that $\varphi ( {\cal A}^H_{\mbmss{\xi}})\in {\cal A}^H_{\mbmss{\xi}}$ holds for all lists $\mbm{\xi}\in J_0^n$ and if, moreover, the $(n+m)$-ary predicate $\sigma_{\varphi}$ defined by $\widetilde{\sigma_{\varphi}}=\{(\mbm{\xi}\,.\,\mbm{\eta})\mid \mbm{\xi}\in J_0^n \,\,\&\,\, \mbm{\eta}\in\varphi ( {\cal A}^H_{\mbmss{\xi}})\}$ is in $J_{n+m}$, then we have $\S\models_{f\langle{S\atop\sigma_{\varphi}} \rangle}\forall \mbm{x} \exists D ( \forall\mbm{y}(D \mbm{y} \leftrightarrow S\mbm{x}\mbm{y}) \land H(\mbm{x}, D )) $. 
\end{enumerate}
\end{lemma}

We have the following for the corresponding finite families $\subseteq J_0^n\times {\cal P}(J_m)$. A part of the proof of Lemma \ref{endlSigma} can be done analogously to the proof of a lemma as described by Fraenkel in \cite{FA37}. 

\begin{lemma}[$choice_h^{n,m}(H)$ for finite families $({\cal A}^H_{\mbmss{\xi}})_{\mbmss{\xi}\in \alpha}$]\label{endlSigma} For any $n\geq 1$, any $m\geq 1$, any Henkin-Asser structure $\S=(J_n)_{n\geq 0}$, any finite $\alpha \in J_n$, any $f$ in ${\rm assgn}(\S)$, and $H(\mbm{x},D)\in {\cal L} ^{(2)}_{\mbmss{x},D}$, 
\[\S\models_{f\langle {A\atop \alpha}\rangle}\forall \mbm{x} \exists D (A\mbm{x}\to H(\mbm{x},D))\] implies 
\[\S\models_{f\langle {A\atop \alpha}\rangle}\exists S \forall \mbm{x} \exists D ( A\mbm{x}\to \forall\mbm{y}(D \mbm{y} \leftrightarrow S\mbm{x}\mbm{y}) \land H(\mbm{x}, D )) .\] 
\end{lemma}
\subsection{Finite supports and adequate decompositions}

The following lemma is a consequence of Equation (1) in Proposition 2.21 (in Part I). The group $\G(P)$ considered in Lemma \ref{EndlichesPfuerH} corresponds to $\G_2$ in Proposition 2.21 if $r\geq 1$ and $D=A_{j_1}^{n_1}$. 

\begin{lemma}[Finite supports and stabilizers for formulas]\label{EndlichesPfuerH} Let $n,m\geq 1$ and $H(\mbm{x},D)$ be in ${\cal L}^{(2)}_{\mbmss{x},D}$. Let $\S$ be any structure of the form $ \S(I,{\G},\I_{\sf 0}^I)$ for some subgroup $\G\subseteq {\G}_{\sf 1}^I$ and let $f$ be in ${\rm assgn}(\S)$. Then, there is a finite set $P\subseteq I$ such that \[\S_{f\langle{\mbmty{x}\atop \mbmty{\xi}}{D\atop\delta}\rangle}(H(\mbm{x},D))=\S_{f\langle{\mbmty{x}\atop \pi(\mbmty{\xi})}{D\atop\delta^\pi}\rangle}(H(\mbm{x},D))\] for all $\pi\in \G(P)$, for all predicates $\delta\in J_m(I,{\G},\I_{\sf 0}^I)$, and for all $\mbm{\xi} \in I^n$.
\end{lemma}

Since we want to decompose the sets $I^n$ --- in particular for $I=\bbbn$ --- by means of $n$-ary predicates of the basic structure $\S_0$, we will say that, for predicates $\alpha, \alpha_1, \ldots, \alpha_l \in {\rm pred}_n(I)$ ($n,l\geq 1$), the set $\{\alpha_1, \ldots, \alpha_l\}$ is a decomposition of $\alpha$ if $\{\widetilde\alpha_1, \ldots, \widetilde{\alpha_{l\,}}\}$ is a decomposition of $\widetilde{\alpha}$. 

\begin{agree}[Notations for decomposing $I^n$] We will introduce some no\-ta\-tions that we will use in the following unless other definitions are given. An overview is given in Table \ref{overview_nota}.
\end{agree}

\begin{table}[t]
\newcounter{liass}

\noindent\framebox{$\!\!\!\!$\parbox{12.05cm}{\begin{list}{$\!\!\!\!\!\!$}
{\usecounter{liass}}
\item\vspace*{-3pt}$\S$ \,\,is a predicate structure $(J_n)_{n\geq 1}$ in ${\sf struc}_{\rm pred}^{({\rm m})}(I)$. 
\item\vspace*{-3pt}$P$ \,\,is a finite subset of $I$.
\item\vspace*{-3pt} $q \,\, =|P|$.
\item\vspace*{-3pt} $\mbm{e} \,\, \in \{1,\ldots,q+1\}^n $.
\item\vspace*{-3pt} ${\rm idx}_{\mbmss{e}} \,\, =\{j\mid 1\leq j\leq n\,\,\&\,\,e_j=q+1\}$.
\item\vspace*{-3pt} ${\sf K}_{{\rm idx}_ {\mbmss{e}}}$ \,\,contains all partitions of ${\rm idx}_ {\mbmss{e}}$.
\item\vspace*{-3pt} ${\cal K} \,\, \in {\sf K}_{ {\rm idx}_ {\mbmty{e}}}$. 
\item\vspace*{-3pt} $\vec {\cal K} \,\, =( K_1, \ldots, K_l)$ if ${\cal K}=\{ K_1, \ldots, K_l\}$ and $\min K_i<\min K_j$ for $i<j$.
\item\vspace*{-3pt} $F_{\cal K}(\mbm{x}) \,\, =\bigwedge_{1\leq v\leq l \atop i\in K_v,\, \min K_v<i}x_{\min K_v}\!=\!x_i \,\,\land\,\, \bigwedge_{1\leq v<w\leq l} \neg (x_{\min K_v}\!=\!x_{\min K_w})$.
\item\vspace*{-3pt} $\alpha_{I,{\cal K}}(\mbm{\xi}) \,\, = \S_{\langle{\mbmty{x} \atop \mbmty{\xi}} \rangle} (F_{\cal K}(\mbm{x}))$.
\item\vspace*{-3pt} $\{\beta_1,\ldots, \beta_{q+1}\}$\,\,is the $P$-adequate partition of $I$ with $\beta_{q+1}=I\setminus P$.
\item\vspace*{-3pt} $F_{\mbmss{e}, {\cal K}}(\mbm{x}) \,\, = B_{e_1}x_1 \land \cdots\land B_{e_n}x_n\land F_{\cal K}(\mbm{x})$. 
\item\vspace*{-3pt} $\alpha_{\mbmss{e}, {\cal K}}^{I,P}(\mbm{\xi}) \,\, = \S_{\langle{\mbmty{x} \atop \mbmty{\xi}} { B_1 \atop\beta_1} {\cdots\atop\cdots} { B_{q+1}\atop\beta_{q+1}}\rangle} (F_{\mbmss{e}, {\cal K}}(\mbm{x}))$ if $\beta_1,\ldots,\beta_{q+1}\in J_1$.
\item\vspace*{-3pt} $\mbm{\mu}$\,\, is a tuple $(\mu_1, \ldots,\mu_n)\in (\beta_{q+1})^n$ with $\mu_i\not=\mu_j$ for $i<j$.
\item\vspace*{-3pt} $P_{\mbmss{\mu}} \,\, =\{\mu_1,\ldots,\mu_n\}$.
\item\vspace*{-3pt} $\mbm{\xi}_{\mbmss{\mu},\mbmss{e}, {\cal K}}^{I,P}$\,\, (shortly, $\mbm{\xi}_{\mbmss{e}, {\cal K}}$) is the tuple
$ (\xi_1, \ldots, \xi_n)$ in $ I^n$ with 

\qquad\quad $\xi_j=\left\{\begin{array}{ll} 
 \nu_{e_j} & \mbox{\rm if $e_j\leq q $ and $\beta_{e_j}=\{ \nu_{e_j} \}$}, \\ \mu_i & \mbox{\rm if ${\cal K}\not=\emptyset$ and $\vec {\cal K}=(K_1,\ldots,K_l)$ and $j \in K_i$}.\\
\end{array}\right.$
\item\label{ASS11}\vspace*{-3pt}
$\alpha^{I,P}_{\mbmss{\mu}}(\mbm{\xi}) \,\, =true$ iff $\mbm{\xi}\in \{\mbm{\xi}_{\mbmss{\mu},\mbmss{e}, {\cal K}}^{I,P} \mid \mbm{e}\in \{1,\ldots,q+1\}^n \,\,\&\,\, {\cal K}\in {\sf K}_{{\rm idx} _{\mbmty{e}}}\}$.
\end{list}}}
\caption{Partitions of $I^n$ and some special notations}\label{overview_nota}
\end{table}

The finite individual supports allow to decompose an individual domain $I$ as follows.
Let $P$ be any finite subset of $ I$. If $P$ is not empty, then let $P=\{\nu_1, \ldots, \nu_q\}$ for some individuals $\nu_1, \ldots, \nu_q\in I$ and $q=|P|$, otherwise, let $q=0$.
Thus, $|P|=q$. For decomposing the domain $I$, we introduce $q+1$ unary predicates. Let, for all $j\in \{1,\ldots,q\}$,
\[\beta_j=\{\nu_j\} \quad\mbox{ and }\quad \beta_{q+1}= I\setminus P.\]
The finite set $\{ \beta_1,\ldots, \beta_{q+1} \}$ is a set of mutually exclusive predicates and thus a partition of $I$ which we call the {\em $P$-adequate partition of the individual domain $I$}.

\begin{lemma}[$P$-adequate partition of the individual domain] For any \linebreak non-empty set $I$ and any finite subset $P$ of $ I$, each predicate $\beta$ in the $P$-adequate partition of $I$ belongs to each Henkin-Asser structure in ${\sf struc}_{\rm pred}^{({\rm m})}(I)$. Moreover, there holds ${\rm sym}_\G(\beta)\supseteq \G(P)$ for each subgroup $\G$ of $\G^I_1$.
\end{lemma}

For any $n\geq 1$ and $q\geq 0$, let $\mbm{e}$ and $\mbm{e}'$ be any $n$-ary tuples of indices in $ \{1,\ldots,q+1\}$. $\mbm{e}$ stands for $(e_1,\ldots,e_n)$ and $\mbm{e}'$ stands for $(e_1',\ldots,e_n')$. 

\setcounter{equation}{0}

\begin{lemma}[$P$-adequately decomposing the set $I^n$]\label{dieBs}$\!$ Let $\S$ be a structure in ${\sf struc}_{\rm pred}^{({\rm m})}(I)$, $f$ be in ${\rm assgn}(\S)$, $n\geq 1$, $P$ be a finite subset of $I$, and $\{ \beta_1,\ldots, \beta_{q+1} \}$ be a $P$-adequate partition of $I$. Let $\mbm{\xi}\in I^n$ and $\pi\in {\G}_{\sf 1}^I(P)$. Moreover, let us assume that $\beta_1,\ldots, \beta_{q+1}$ belong to $\S$ and that $\bar f$ is the assignment $f\langle{ B_1 \atop\beta_1} {\cdots\atop \cdots} { B_{q+1}\atop\beta_{q+1}}\rangle$. 
Then, the statements {\rm (\ref{dieBs2})} and {\rm (\ref{dieBs20})} and,  for each $\mbm{e}\in \{1,\ldots,q+1\}^n$, the statements {\rm (\ref{dieBs2b})} and {\rm (\ref{dieBs0})} hold. 

\begin{equation}\label{dieBs2}\S\models_{\bar f}
\bigvee_{\mbmss{e}\in \{1,\ldots,q+1\}^n} (B_{e_1}x_1 \land \cdots\land B_{e_n}x_n )\end{equation}

\begin{equation}\label{dieBs20}\S\models_{\bar f}\bigwedge_{\mbmss{e}, \mbmss{e'}\in \{1,\ldots,q+1\}^n \atop \mbmss{e}\not= \mbmss{e}'} (B_{e_1}x_1 \land \cdots\land B_{e_n}x_n \to \neg (B_{e_1'}x_1 \land \cdots\land B_{e_n'}x_n ))\end{equation}

\begin{equation}\label{dieBs2b}\S_{\bar f\langle{\mbmty{x}\atop\mbmty{\xi}} \rangle}
(B_{e_1}x_1 \land \cdots\land B_{e_n}x_n )=\S_{\bar f\langle{\mbmty{x}\atop\pi(\mbmty{\xi} )} \rangle}
(B_{e_1}x_1 \land \cdots\land B_{e_n}x_n )\end{equation}

\vspace{0.1cm}
 If $\S$ is a Henkin-Asser structure, then the predicates $\beta_1,\ldots, \beta_{q+1}$ and consequently the predicate $\beta_{e_1} \times \cdots\times \beta_{e_n}$ given by 
  \begin{equation}\label{dieBs0}\beta_{e_1} \times \cdots\times \beta_{e_n}=\alpha_{\S,B_{e_1}x_1 \land \cdots\land B_{e_n}x_n,\mbmss{x},\bar f}\end{equation} belong to $\S$. 
\end{lemma}

By (\ref{dieBs2b}) in Lemma \ref{dieBs}, we have ${\rm sym}_\G(\beta_{e_1} \times \cdots\times \beta_{e_n})\supseteq \G(P)$ for each subgroup $\G$ of $\G^I_1$. For each $P$-adequate partition $\{ \beta_1,\ldots,\beta_{q+1} \}$ of $I$, by the statements (\ref{dieBs2}) and (\ref{dieBs20}) in Lemma \ref{dieBs}, the set $\{\beta_{e_1} \times \cdots\times \beta_{e_n} \mid \mbm{e}\in \{1,\ldots,q+1\}^n \}$ is a finite partition of $I^n$ that we call the {\em $P$-adequate partition of $I^n$}. For the proof of $\S_0\models choice_h(H)$ for $H$ in ${\cal L}^{(2)}_{\mbmss{x},D}$, it is useful to refine each $P$-adequate partition of each product $I^n$ (if $n\geq 2$) since the fact that two or more components of certain $n$-tuples $(\xi_1,\ldots,\xi_n)\in I^n$ can be equal has to be taken into account. Let $\xi_1,\ldots,\xi_n$ be any list of individuals in $I$. The next lemma deals with the following questions. What happens when we change the values $\xi_1,\ldots,\xi_n$ by a permutation $\pi\in \G_{\sf 1}^I(P)$? What happens when we change some of the values in $I\setminus P$? We know that, for all $i,j\in \{1,\ldots, n\}$, $\xi_{i}=\xi_{j}$ holds iff $\pi(\xi_{i})=\pi(\xi_{j})$ holds. Because of this property, we will consider a partition ${\cal K}$ of the set ${\rm idx}$ containing all indices $j\in \{1, \ldots, n\} $ satisfying $\xi_j\in I\setminus P$. Let $\{\beta_1,\ldots, \beta_{q+1}\}$ be the $P$-adequate partition of $I$ with $\beta_{q+1}=I\setminus P$. If the index set ${\rm idx}=\{j\mid 1\leq j\leq n\,\,\&\,\,\xi_j\in \beta_{q+1}\}$ is not empty, such a partition ${\cal K}$ should satisfy $\bigcup{\cal K}={\rm idx}$ and any two indices $i$ and $j$ in ${\rm idx}$ should be in one and the same set $K\in{\cal K}$ iff $\xi_{i}=\xi_{j}$ and $\xi_{i},\xi_{j}\in \beta_{q+1}$. Later, we will describe such a ${\cal K}$ by exactly one ordered tuple $\vec {\cal K}=(K_1, \ldots, K_l)$ which implies ${\cal K}=\{K_1, \ldots, K_l\}$ ($l\leq |{\rm idx}|$). 
This means that for such a list $\xi_1,\ldots,\xi_n$ and such a ${\cal K}$ in the case that $|{\rm idx}|>1$, the formula $F_ {\cal K}(\mbm{x)}$ uniquely defined by 
\[F_{\cal K}(\mbm{x})=\bigwedge_{1\leq v\leq l}\,\bigwedge_{i\in K_v\atop \min K_v<i}x_{\min K_v}=x_i \,\,\land \bigwedge_{1\leq v<w\leq l}\, \neg (x_{\min K_v}=x_{\min K_w})\] will be true in any structure $\S\in {\sf struc}_{\rm pred}^{({\rm m})}(I)$ if $(\xi_1,\ldots,\xi_n)$ is assigned to $\mbm{x}$. 
Moreover, suitable to the abbreviations given in Section 2 (of Part I), we will use the formulas $F_{\{\{i\}\}}$ and $F_{\emptyset}$ given by
\[F_{\{\{i\}\}}(\mbm{x})=x_i=x_i \quad \mbox{ and }\quad F_{\emptyset}(\mbm{x})=x_1=x_1\] in the trivial cases where $|{\rm idx}|\leq 1$ and the partitions are $\{\{i\}\}$ for some $i\leq n$ and $\emptyset$, respectively. 
For decomposing any non-empty index set ${\rm idx}\subseteq \{1,\ldots,n\}$, we consider, for each $l\in\{1,\ldots,| {\rm idx}|\}$, certain $l$-tuples $(K_1,\ldots,K_l)$ of pairwise disjoint non-empty subsets $K_v\subseteq{\rm idx}$ with $K_1\cup\cdots\cup K_l={\rm idx}$. The corresponding sets $\{K_1,\ldots,K_l\}$ called {\em $l$-partitions of ${\rm idx}$} are collected in ${\sf K}_{{\rm idx},l} $ as follows. Let ${\cal P}^+(\bbbn)$ be the set $\{K\subseteq \bbbn\mid K\not=\emptyset\}$ and \[{\sf K}_{{\rm idx},l}=\{\{K_1,\ldots,K_l\}\subset {\cal P}^+(\bbbn)\mid\]\[\hspace{3cm} \bigcup\limits_{i=1}^lK_i={\rm idx}\,\,\,\&\,\,\, ( 1\leq i<j\leq l \Rightarrow K_i\cap K_j=\emptyset)\}\] be the finite set of all $l$-partitions of ${\rm idx}$. Moreover, let ${\sf K}_{\rm idx}$ be the finite set of all partitions of ${\rm idx}$ given by ${\sf K}_{\rm idx}= \bigcup_{l=1}^{|{\rm idx}|}\,\,{\sf K}_{{\rm idx},l}$. The only {\em partition of $\emptyset$} is the empty set. Thus, let ${\sf K}_{\emptyset}=\{\emptyset\}$. Let us remark that, in case of ${\rm idx}= \{1,\ldots,n\}$, the number $|{\sf K}_{\rm idx}|$ of sets in ${\sf K}_{\rm idx}$ is the Bell number $B_n$ (for the definition see, e.g., \cite{DudReMa}). For any $l$-partition ${\cal K}$ of a subset of $\{1,\ldots,n\}$, let $\vec {\cal K}$ be the $l$-tuple $( K_1, \ldots, K_l)$ such that ${\cal K}=\{ K_1, \ldots, K_l\}$ holds and $\min K_i<\min K_j$ is satisfied for all $i$ and $j$ with $1\leq i<j\leq l$. For any non-empty set $I$ and any partition ${\cal K}$ of any subset ${\rm idx}\subseteq \{1,\ldots,n\}$, let the predicate $ \alpha_{I,{\cal K}}\in {\rm pred}_n(I)$ be defined by 
\[\alpha_{I,{\cal K}}(\mbm{\xi}) = \S_{\langle{\mbmty{x} \atop \mbmty{\xi}} \rangle} (F_{\cal K}(\mbm{x}))\] 
for any $\mbm{\xi}\in I^n$ where $\S$ is any structure in ${\sf struc}_{\rm pred}^{({\rm m})}(I)$. 

\setcounter{equation}{0}
\begin{lemma}[${\rm id}$-and-${\rm idx}$ adequately decomposing $I^n$]\label{Zerl_idx} Let $n\geq 1$, $\S$ be any structure in ${\sf struc}_{\rm pred}^{({\rm m})}(I)$, $f\in {\rm assgn}(\S)$, and ${\rm idx}$ be a non-empty subset of $ \{1,\ldots,n\}$. Let $\{\alpha_1,\ldots, \alpha_{ |{\sf K}_{ {\rm idx}}|}\}$ be the set $\{\alpha_{I,{\cal K}}\mid {\cal K} \in {\sf K}_{{\rm idx}}\}$ and $\mbm{\xi}\in I^n$. If $\alpha_s$ belongs to $\S$ for some $s\in \{1,\ldots, |{\sf K}_{ {\rm idx}}|\}$, then, for the  assignment $f'=f\langle{A_s\atop \alpha_{s}}\rangle$ and any $\pi\in {\G}_{\sf 1}^I$, we have $ \S_{f'\langle{\mbmty{x}\atop\mbmty{\xi}} \rangle}(A_s\mbm{x})=\S_{f'\langle{\mbmty{x}\atop\pi(\mbmty{\xi} )} \rangle}(A_s\mbm{x})$ and, for any ${\cal K} \in {\sf K}_{{\rm idx}}$, the statement $(\ref{dieGs4})$. If $\S$ moreover is a Henkin-Asser structure, then, for any assignment $f'=f\langle{A_1\atop \alpha_{1}}{\cdots\atop\cdots}{A_{|{\sf K}_{ {\rm idx}}|} \atop \alpha_{{|{\sf K}_{ {\rm idx}}|}}} \rangle$, $\S\models_{f'} \bigvee_{s=1,\ldots, |{\sf K}_{ {\rm idx}}|} A_s\mbm{x}$ and $(\ref{dieGs4a})$ are satisfied.
\begin{equation}\label{dieGs4}\S_{f\langle{\mbmty{x}\atop \mbmty{\xi}}\rangle} ( F_{\cal K}(\mbm{x}))=\S_{f\langle{\mbmty{x}\atop \pi(\mbmty{\xi})}\rangle} ( F_{\cal K}(\mbm{x})) \qquad ( \mbox{for any ${\cal K} \in {\sf K}_{{\rm idx}}$})
\end{equation} 
\begin{equation}\label{dieGs4a}\S\models_f\bigvee_ {{\cal K} \in {\sf K}_{{\rm idx}}} F_{\cal K}(\mbm{x}) 
\end{equation}\end{lemma}
\begin{corollary} For each ${\rm idx}\subseteq \{1,\ldots,n\}$, $\{\alpha_{I,{\cal K}}\mid {\cal K} \in {\sf K}_{{\rm idx}}\}$ is a finite partition of $I^n$ that is  invariant with respect to the permutations in ${\G}_{\sf 1}^I$. Moreover, for each $\alpha\subseteq I^n$, $\{\alpha_{I,{\cal K}}\cap \alpha \mid {\cal K} \in {\sf K}_{{\rm idx}}\}$ is a finite decomposition of $\alpha$.\end{corollary} 

For any ${\rm idx}\subseteq \{1,\ldots,n\}$, the partition $\{\alpha_{I,{\cal K}}\mid {\cal K} \in {\sf K}_{{\rm idx}}\}$ is called the {\em ${\rm id}_{\rm idx}$-adequate partition of $I^n$}. For ${\rm idx}=\{1,\ldots,n\}$ this partition is called the {\em ${\rm id}$-adequate partition of $I^n$}. If $\{\alpha_{I,{\cal K}}\cap \alpha \mid {\cal K} \in {\sf K}_{{\rm idx}}\}$ is a partition of $\alpha\subseteq I^n$, then it is called the {\em ${\rm id}_{\rm idx}$-adequate partition of $\alpha$}.

Now, let $P$ be any finite subset of a set $I$, $q =|P|$, $\{\beta_1,\ldots, \beta_{q+1}\}$ be the $P$-adequate partition of $I$ with $\beta_{q+1}=I\setminus P$, and let $\mbm{e}\in \{1,\ldots,q+1\}^n $. Then, let ${\rm idx}_{\mbmss{e}}=\{j\mid 1\leq j\leq n\,\,\&\,\,e_j=q+1\}$ and $F_{\mbmss{e}, {\cal K}}(\mbm{x})$ be the formula given by \[F_{\mbmss{e}, {\cal K}}(\mbm{x})=B_{e_1}x_1 \land \cdots\land B_{e_n}x_n\land F_{\cal K}(\mbm{x})\] for all ${\cal K} \in {\sf K}_{ {\rm idx}_ {\mbmty{e}}}$. Let $\S$ be a Henkin-Asser structure in ${\sf struc}_{\rm pred}^{({\rm m})}(I)$, $f\in{\rm assgn}(\S)$, and $\bar f=f\langle{ B_1 \atop\beta_1} {\cdots\atop \cdots} { B_{q+1}\atop\beta_{q+1}}\rangle$. Consequently, $\{\{\mbm{\xi}\in I^n\mid \S\models _{\bar f\langle{\mbmty{x}\atop\mbmty{\xi}}\rangle}F_{\mbmss{e}, {\cal K}}(\mbm{x})\}\mid{\cal K} \in {\sf K}_{ {\rm idx}_ {\mbmty{e}}}\}$ is the ${\rm id}_{{\rm idx}_{\mbmty{e}}}$-adequate partition of $\beta_{e_1}\times \cdots\times\beta_{e_n}$. Thus, we can combine the ${\rm id}_{{\rm idx}_{\mbmty{e}}}$-adequate partition of a subset $\beta_{e_1}\times \cdots\times\beta_{e_n}$ with the $P$-adequate partition of $I^n$ given by $\{\beta_{e_1'} \times \cdots\times \beta_{e_n'} \mid \mbm{e}'\in \{1,\ldots,q+1\}^n \}$ in order to get a refinement of the $P$-adequate partition of $I^n$. This means that we take the partitions $\{ \alpha_{1,1},\ldots, \alpha_{1,{l_1}} \}= \{\beta_{e_1'} \times \cdots\times \beta_{e_n'} \mid \mbm{e}'\in \{1,\ldots,q+1\}^n \}$ and $\{ \alpha_{2,1},\ldots, \alpha_{2,{l_2}} \}=\{\alpha_{I,{\cal K}} \mid {\cal K} \in {\sf K}_{{\rm idx}_{\mbmss{e}}}\} $ and we get the ${\rm id}_{{\rm idx}_{\mbmty{e}}}$-adequate partition $\{\alpha_{I,{\cal K}}\cap(\beta_{e_1}\times \cdots\times\beta_{e_n}) \mid {\cal K} \in {\sf K}_{{\rm idx}_{\mbmty{e}}}\}$ for $\alpha_{1,i}= \beta_{e_1}\times \cdots\times\beta_{e_n}$ since it does not contain $\{\emptyset\} $. More generally, for all $ \beta_{e_1}\times \cdots\times\beta_{e_n}$, we need the ${\rm id}_{{\rm idx}_{\mbmty{e}}}$-adequate partition. For all $n\geq 1$ and, in particular, for each $n\geq 2$, we want to refine the $P$-adequate partition of $I^n$ and consider --- as far as necessary --- common refinements of the $P$-adequate partition of $I^n$ combined with the ${\rm id}_{{\rm idx}_{\mbmty{e}}}$-adequate partitions of all possible sets $\beta_{e_1}\times \cdots\times\beta_{e_n}$ in the $P$-adequate partition where $\mbm{e}\in \{1,\ldots,q+1\}^n$. Each partition suitable for the following construction is the result of refining the $P$-adequate partition of $I^n$ and all resulting partitions by means of the ${\rm id}_{{\rm idx}_{\mbmty{e}}}$-adequate partition of a set $\beta_{e_1}\times \cdots\times\beta_{e_n}$ stepwise for each $\mbm{e}\in \{1,\ldots,q+1\}^n$. Let 
\[\alpha_{\mbmss{e}, {\cal K}}^{I,P}(\mbm{\xi}) = \S_{\bar f\langle{\mbmty{x} \atop \mbmty{\xi}} \rangle} (F_{\mbmss{e}, {\cal K}}(\mbm{x})).\] 
 
Since, by Lemma \ref{dieBs}, $\S_{\bar f\langle{\mbmty{x} \atop \mbmty{\xi}} \rangle} ( \bigwedge_{\mbmss{e}, \mbmss{e'}\in \{1,\ldots,q+1\}^n \atop \mbmss{e}\not= \mbmss{e}'} (F_{\mbmss{e}, {\cal K}}(\mbm{x}) \to \neg (B_{e_1'}x_1 \land \cdots\land B_{e_n'}x_n )))$ is $true$ for any Henkin-Asser structure $\S\in{\sf struc}_{\rm pred}^{({\rm m})}(I)$, \[\{\alpha_{\mbmss{e}, {\cal K}}^{I,P} \mid \mbm{e}\in \{1,\ldots,q+1\}^n \,\,\&\,\, {\cal K}\in {\sf K}_{{\rm idx} _{\mbmty{e}}}\}\] is a finite partition of $I^n$ that we call the {\em $P$-and-${\rm id}$ adequate partition of $I^n$}.

As a consequence of Lemma \ref{Zerl_idx} we get the following description.

\setcounter{equation}{0}
\begin{lemma}[$P$-and-${\rm id}$ adequate partition and stability]\label{Lemma_P_id} 
Let $\S$ be any \linebreak Henkin-Asser structure in ${\sf struc}_{\rm pred}^{({\rm m})}(I)$ and $f$ be any assignment in ${\rm assgn}(\S)$.
Let $P$ be a finite subset of $I$ and let $\{\beta_1,\ldots, \beta_{q+1}\}$ be the $P$-adequate partition of $I$ with $\beta_{q+1}=I\setminus P$. Let $n\geq 1$ and $\{\alpha_1,\ldots, \alpha_k\}$ be the $P$-and-${\rm id}$ adequate partition of $I^n$. Moreover, let $f'=f\langle{A_1\atop \alpha_1}{\cdots\atop\cdots}{A_k \atop \alpha_k} \rangle$, $\bar f=f\langle{ B_1 \atop\beta_1} {\cdots\atop\cdots} { B_{q+1}\atop\beta_{q+1}}\rangle$, and $\pi\in {\G}_{\sf 1}^I(P)$. Then, we have
\[\S \models_{f'}\bigvee_{s=1,\ldots, k} A_s\mbm{x} \quad \mbox{ and } \quad \S_{f'\langle{\mbmty{x}\atop\mbmty{\xi}}\rangle}(A_s\mbm{x})
=\S_{f'\langle{\mbmty{x}\atop\pi(\mbmty{\xi} )}\rangle} (A_s\mbm{x})\]
for all $s\in \{1,\ldots, k\}$ and thus the statement $(\ref{dieGs4Z})$ for each $\mbm{e}\in \{1,\ldots,q+1\}^n $ and each 
 ${\cal K} \in {\sf K}_{ {\rm idx}_ {\mbmty{e}}}$ and $(\ref{Zerlegung})$.
\begin{equation}\label{dieGs4Z}\S_{\bar f\langle{\mbmty{x}\atop\mbmty{\xi}}\rangle} ( F_{\mbmss{e}, {\cal K}}(\mbm{x}))=\S_{\bar f\langle{\mbmty{x}\atop\pi(\mbmty{\xi})}\rangle} ( F_{\mbmss{e}, {\cal K}}(\mbm{x})) 
\end{equation} 

\begin{equation}\label{Zerlegung}
\S\models_{\bar f} \bigvee_{\mbmss{e}\in \{1,\ldots,q+1\}^n}\, \bigvee_{{\cal K} \in {\sf K}_{ {\rm idx}_ {\mbmty{e}}}} F_{\mbmss{e}, {\cal K}}(\mbm{x})
\end{equation} 
\end{lemma}

For defining a choice set for these partitions, let $\mbm{\mu}$ be any tuple $(\mu_1, \ldots,\mu_n)$ of individuals in $I \setminus P$ satisfying $\mu_i\not=\mu_j$ for $i\leq i<j\leq n$ and $l= |{\cal K}|$. $\mbm{\xi}_{\mbmss{\mu},\mbmss{e}, {\cal K}}^{I,P}$ denotes the tuple $(\xi_1, \ldots, \xi_n)$ in $\alpha_{\mbmss{e}, {\cal K}}^{I,P}$ determined as follows where $(K_1,\ldots,K_l)$ is the ordered tuple $\vec {\cal K}$ in case that $l\in \{1,\ldots,n\}$. 
\[\xi_j=\left\{\begin{array}{ll} \nu_{e_j} & \mbox{\rm if $e_j\leq q $ and $\beta_{e_j}=\{ \nu_{e_j} \}$ \quad (and thus $j \not\in \bigcup_{i=1}^l K_i$)}, \\ \mu_i & \mbox{\rm if $l\geq 1 $ and $j \in K_i$ for some $i\leq l$}.\\
\end{array}\right.\] Moreover, let $\alpha^{I,P}_{\mbmss{\mu}}$ be the predicate in ${\rm pred}_n(I)$ satisfying \[\widetilde {\alpha^{I,P}_{\mbmss{\mu}}}=\{\mbm{\xi}_{\mbmss{\mu},\mbmss{e}, {\cal K}}^{I,P} \mid \mbm{e}\in \{1,\ldots,q+1\}^n \,\,\&\,\, {\cal K}\in {\sf K}_{{\rm idx} _{\mbmty{e}}}\}\] and let $P_{\mbmss{\mu}}=\{\mu_1,\ldots,\mu_n\}$.

\begin{lemma}[A choice set for a $P$-and-${\rm id}$ adequate partition]\label{L_P_id} 
 Let $n\geq 1$ and $\mbm{\mu}$ be any tuple $(\mu_1, \ldots,\mu_n)$ of individuals in $I \setminus P$ with $\mu_i\not=\mu_j$ for $i\leq i<j\leq n$. Then, $\widetilde{\alpha^{I,P}_{\mbmss{\mu}}}$ is a choice set for the $P$-and-${\rm id}$ adequate partition of $I^n$.
For all subgroups $\G\subseteq \G_{\sf 1}^I$, we have 
\[\quad{\rm sym}_{\G}(\{\mbm{\xi}_{\mbmss{\mu},\mbmss{e}, {\cal K}}^{I,P}\})\supseteq \G(P\cup P_{\mbmss{\mu}})\] for any $ \mbm{e}\in \{1,\ldots,q+1\}^n$ and ${\cal K}\in {\sf K}_{{\rm idx} _{\mbmty{e}}}$ and thus 
\[{\rm sym}_{\G}(\alpha_{\mbmss{\mu}}^{I,P} )\supseteq \G(P\cup P_{\mbmss{\mu}}).\] 
\end{lemma}

\subsection{Describing consequences of swapping two tuples of individuals}\label{DefinPerm}

Let $P$ be a finite subset of a domain $I$ and $q=|P|$. Let $n\geq 1$, $\mbm{e} \in \{1,\ldots,q+1\}^n$, ${\cal K} \in {\sf K}_{{\rm idx}_{\mbmss{e}}}$, and $\mbm{\mu}$ be a tuple $(\mu_1,\ldots,\mu_n)$ in $(I\setminus P)^n$ whose components are pairwise different. For ${\cal K}\not= \emptyset$, let $\vec {\cal K}=(K_1,\ldots,K_{|{\cal K}|})$, and $M_{\mbmss{e}, {\cal K}} = \{\mu_1, \ldots,\mu_{|{\cal K}|}\}$. 

Then, the tuple $\mbm{\xi}_{\mbmss{\mu},\mbmss{e}, {\cal K}}^{I,P}$ is uniquely determined. Here, it will be denoted by $\mbm{\xi}_{\mbmss{e}, {\cal K}}$ and $(\xi^{\mbmss{e}, {\cal K}}_1,\ldots,\xi^{\mbmss{e}, {\cal K}}_n)$. If $\mbm{\xi}_{\mbmss{e}, {\cal K}}\not \in P^n$, then ${\cal K}\not= \emptyset$ and for certain indices $i_1,\ldots,i_{|{\cal K}|}$ in ${\rm idx}_{\mbmss{e}}$, we have $ \xi_{i_1}^{\mbmss{e}, {\cal K}}=\mu_1, \ldots,\xi_{i_{|{\cal K}|}}^{\mbmss{e}, {\cal K}}=\mu_{|{\cal K}|}$ and for $i\not\in {\rm idx}_{\mbmss{e}}$ we have $\xi_{i}^{\mbmss{e}, {\cal K}}\in P$. 

Let $\{\beta_1,\ldots, \beta_{q+1}\}$ be the $P$-adequate partition of $I$ with $\beta_{q+1}=I\setminus P$ and $\mbm{\xi}$ be any $n$-tuple $(\xi_1,\ldots,\xi_n)\in I^n$ satisfying the relationship $\S\models_{f\langle{ B_1 \atop\beta_1} {\cdots\atop\cdots} { B_{q+1}\atop\beta_{q+1}}{\mbmty{x}\atop \mbmty{\xi}}\rangle}F_{\mbmss{e}, {\cal K}}(\mbm{x})$. If $\mbm{\xi}_{\mbmss{e}, {\cal K}} \in P^n$, then $\mbm{\xi} \in P^n$ and $\mbm{\xi}_{\mbmss{e}, {\cal K}} =\mbm{\xi}$ by definition.

If we assume ${\cal K}\not=\emptyset$ and we compare $\mbm{\xi}$ with $\mbm{\xi}_{\mbmss{e}, {\cal K}}$, then $\{\xi_{i_1},\ldots, \xi_{i_{|{\cal K}|}}\}$ is the counterpart to $\{\xi_{i_1}^{\mbmss{e}, {\cal K}},\ldots, \xi_{i_{|{\cal K}|}}^{\mbmss{e}, {\cal K}}\}$ and we have $\{\xi_i\mid i\in {\rm idx}_{\mbmss{e}}\}=\{\xi_{i_1},\ldots, \xi_{i_{|{\cal K}|}}\}$ and $\{\xi_i^{\mbmss{e}, {\cal K}}\mid i\in {\rm idx}_{\mbmss{e}}\}=\{\xi_{i_1}^{\mbmss{e}, {\cal K}},\ldots, \xi_{i_{|{\cal K}|}}^{\mbmss{e}, {\cal K}}\}=M_{\mbmss{e}, {\cal K}}$ and the two sets need not be equal.
On the other hand, the two sets need not be disjoint and it is a question what effects the application of a permutation $\pi$ has if $\pi$ assigns $\mbm{\xi}$ to $\mbm{\xi}_{\mbmss{e}, {\cal K}}$. For ${\cal K}\not= \emptyset$, we want to recursively define a list $\zeta_1,\ldots,\zeta_{s_{\mbmty{\xi}}}$ of pairwise different individuals satisfying $\{\zeta_1,\ldots,\zeta_{s_{\mbmty{\xi}}}\}=M_{\mbmss{e}, {\cal K}}\cup \{\xi_i\mid i\in {\rm idx}_{\mbmss{e}}\}$ such that it is possible to define a $\pi_{\mbmss{\xi}}\in \G^I_{\sf 1}(I\setminus \{\zeta_1,\ldots,\zeta_{s_{\mbmty{\xi}}}\})$ with $\mbm{\xi}=\pi_{\mbmss{\xi}} (\mbm{\xi}_{\mbmss{e}, {\cal K}} )$. For ${\cal K}= \emptyset$, let $\pi_{\mbmss{\xi}}$ be in $\G^I_{\sf 1}(I)$ and $s_{\mbmss{\xi}}=0$. In the latter case, $\pi_{\mbmss{\xi}}$ is the only mapping in $\G^I_{\sf 1}(I)$, it is the identity bijection.
Important for us will be the fact that the permutation $\pi_{\mbmss{\xi}}$ defined as follows is in $\G^I_{\sf 1}(P)$ and that its application to any $m$-tuple of individuals given by $\pi_{\mbmss{\xi}}(\mbm{\eta}_0)=\mbm{\eta}$ can be uniquely described by a finite formula in ${\cal L}^{(2)}_{\mbmss{x}_0,\mbmss{y}_0, \mbmss{x},\mbmss{y}}$ under the assignment $\langle {\mbmss{x}_0\atop\mbmss{\xi}_{\mbmty{e}, {\cal K}}}{\mbmss{y}_0\atop\mbmss{\eta}_0}{\mbmss{x}\atop\mbmss{\xi}}{\mbmss{y}\atop\mbmss{\eta}}\rangle $.
For ${\cal K}\not= \emptyset$, let 
\[\!\!\!\!\!\!\!\!\zeta_j= \left\{\begin{array}{ll} \mu_j &\mbox{ if }1\leq j \mbox{ and } j \leq |{\cal K}|, \\
\xi_t &\mbox{ if }0<|{\cal K}|<j \leq s_{\mbmss{\xi}} \mbox{ and } t= \min\{i\in {\rm idx}_{\mbmss{e}}\mid \xi_i\not\in \{\zeta_1,\ldots,\zeta_{j-1}\}\}.
\end{array}\right.\]
Let $\pi_{\mbmss{\xi}}\in \G^I_{\sf 1}(P)$ be a permutation defined step-by-step by 
\[\pi_{\mbmss{\xi}}(\zeta)=
\left\{\begin{array}{ll} \zeta&\mbox{ if } \zeta\in I \mbox{ and } \zeta\not = \zeta_j \mbox{ for all } j\leq s_{\mbmss{\xi}},\\
\xi_{i_j}&\mbox{ if } \zeta=\zeta_j \mbox{ for some } j \mbox{ with }1\leq j \leq |{\cal K}|, \\
\zeta_{t} &\mbox{ if } \zeta = \zeta_j \mbox{ for some $j$ with }0<|{\cal K}|< j \leq s_{\mbmss{\xi}} \mbox{ and }\\& \hspace{0.48cm} t=\min\{i \in\{1,\ldots,|{\cal K}|\}\mid \zeta_i\not\in \{\pi_{\mbmss{\xi}}(\zeta_1),\ldots,\pi_{\mbmss{\xi}}(\zeta_{j-1})\}\}. \end{array}\right.\]
Thus, we have
 $\mbm{\xi}=\pi_{\mbmss{\xi}} (\mbm{\xi}_{\mbmss{e}, {\cal K}} )$ and $\mbm{\xi}_{\mbmss{e}, {\cal K}} =\pi_{\mbmss{\xi}}^{-1}( \mbm{\xi})$. 
If $\{ \xi_1^{\mbmss{e}, {\cal K}}, \ldots, \xi_n^{\mbmss{e}, {\cal K}},\xi_1,\ldots,\xi_n\}\setminus P$ is not empty, then $|{\cal K}|>0$. If $|{\cal K}|>0$ and $s_{\mbmss{\xi}}= |{\cal K}|$ hold, then the set $\{\zeta_1,\ldots,\zeta_{s_{\mbmty{\xi}}}\}$ contains exactly the pairwise different individuals $ \xi_{i_1}^{\mbmss{e}, {\cal K}}, \ldots, \xi_{i_{|{\cal K}|}}^{\mbmss{e}, {\cal K}}$. Otherwise, for $s_{\mbmss{\xi}}> |{\cal K}|>0$, it contains $s_{\mbmss{\xi}}$ pairwise different individuals $ \xi_{i_1}^{\mbmss{e}, {\cal K}}, \ldots, \xi_{i_{|{\cal K}|}}^{\mbmss{e}, {\cal K}},\xi_{j_1},\ldots,\xi_{j_{s_{\mbmty{\xi}}-|{\cal K}|}}$ for some indices $j_1, \ldots, j_{s_{\mbmty{\xi}}-|{\cal K}|}\in\{i_1,\ldots,i_{|{\cal K}|}\}$. Consequently, $\pi_{\mbmss{\xi}}$ is completely determined by a permutation of the set of these individuals that can also be described by a permutation of $\{1,\ldots,|{\cal K}|\}\subseteq \{1,\ldots, n\}$ and $\{1,\ldots,|{\cal K}|, |{\cal K}|+1,\ldots, s_{\mbmty{\xi}}\}\subseteq \{1,\ldots, 2n\}$, respectively. This implies that the relationship between the tuples $(\mbm{\xi}_{\mbmss{e}, {\cal K}} \,.\,\mbm{\eta}_0)$ and $(\mbm{\xi} \,.\,\pi_{\mbmss{\xi}}(\mbm{\eta}_0))$ can be described by one formula $swap_{\mbmss{e}}$ --- for fixed $n,m\geq 1$ --- under an assignment of the form $\langle {\mbmss{x}_0\atop\mbmss{\xi}_{\mbmty{e}, {\cal K}}}{\mbmss{y}_0\atop\mbmss{\eta}_0}{\mbmss{x}\atop\mbmss{\xi}}{\mbmss{y}\atop\mbmss{\eta}}\rangle $.
We will use $\mbm{x}_0=(x_1^{(0)},\ldots,x_n^{(0)})$, $\mbm{x}=(x_1,\ldots,x_n)$, $\mbm{y}_0=(y_1^{(0)},\ldots,y_m^{(0)})$, and in particular $\mbm{y}=(y_1,\ldots,y_m)$ for describing $\mbm{\eta}=\pi_{\mbmss{\xi}} (\mbm{\eta}_0)$. 

\vspace{0.2cm}

\noindent For each $\mbm{e}$ with ${\rm idx}_{\mbmss{e}}=\emptyset$, $swap_{\mbmss{e}}(\mbm{x}_0,\mbm{y}_0,\mbm{x},\mbm{y})$ stands for $\mbm{x}= \mbm{x}_0 \land \mbm{y}= \mbm{y}_0 $. 

\vspace{0.2cm}
\noindent For each $\mbm{e}$ with ${\rm idx}_{\mbmss{e}}\not=\emptyset$, $swap_{\mbmss{e}}(\mbm{x}_0,\mbm{y}_0,\mbm{x},\mbm{y})$ stands for
\[(\mbm{x}= \mbm{x}_0 \to \mbm{y}= \mbm{y}_0 ) \,\,\land\,\, (\mbm{x}\not = \mbm{x}_0 \to\hspace*{6.7cm}\]
\[\land \bigwedge_ { k=1}^m
(( \bigwedge_ {i\in {\rm idx}_{\mbmss {e}}} (y^{(0)}_{k} \not=x^{(0)}_i\land y^{(0)}_{k} \not=x_i) \to y_{k}=y^{(0)}_{k})\hspace*{4.0cm}
\]\[\land \bigwedge_ {i\in {\rm idx}_{\mbmss {e}}} (y^{(0)}_k =x^{(0)}_i \to y_{k}=x_i )\]
\[\land \bigwedge_ {i\in {\rm idx}_{\mbmss {e}}}(y^{(0)}_k =x_i \land \bigwedge_ { j=1}^n x_i \not= x_j^{(0)}\]
\[\hfill
\to \bigvee_{1\leq t\leq i} \,\, \bigvee_{i_0\leq n} enum(t, i,\mbm{x},\mbm{x}_0) \land enum(t, i_0,\mbm{x}_0,\mbm{x}) \land y_{k} = x_{i_0}^{(0)})
 )).\]
The formulas $enum(t, i,\mbm{x},\mbm{x}_0)$ and $enum(t, i_0,\mbm{x}_0,\mbm{x})$ belong to ${\cal L}^{(2)}_{\mbmss{x}, \mbmss{x}_0}$. Let us use
$enum(t, i,\mbm{x},\mbm{x}_0)$ to describe that $ \xi_i$ is the $t{\rm th}$ new element in the sequence $(\xi_1, \ldots , \xi_i)$ that does not occur in the list $\mbm{\xi}_{\mbmss{e}, {\cal K}}$ which means that $ \xi_i$ is not in $ \{\xi^{\mbmss{e}, {\cal K}}_1,\ldots,\xi^{\mbmss{e}, {\cal K}}_n,\xi_1, \ldots , \xi_{i-1} \}$ and that $| \{\xi^{\mbmss{e}, {\cal K}}_1,\ldots,\xi^{\mbmss{e}, {\cal K}}_n,\xi_1, \ldots , \xi_{i-1} \}|=|{\cal K}|+t-1$. Let $enum(t, i_0,\mbm{x}_0,\mbm{x})$ say that $ \xi^{\mbmss{e}, {\cal K}}_{i_0}$ is not in $ \{\xi_1,\ldots,\xi_n,\xi_1^{\mbmss{e}, {\cal K}}, \ldots , \xi_{i_0-1} ^{\mbmss{e}, {\cal K}}\}$ and that it is the $t^{\rm th}$ new element in the sequence $(\xi^{\mbmss{e}, {\cal K}}_1,\ldots,\xi^{\mbmss{e}, {\cal K}}_{i_0})$ that does not occur in the list $\mbm{\xi}$. Let $\mbm{v}$ and $\mbm{w}$ stand for the list $(v_1, \ldots , v_n)$ and $(w_1, \ldots , w_n)$, respectively. 
\begin{description}
\item $enum(1, i,\mbm{v},\mbm{w})$ is the formula 
\[\bigwedge_{j=1}^{n} v_i\not= w_j\land \bigwedge_{ s < i}\,\,\,\bigvee_{j=1}^{n} v_{s}= w_j.\hspace*{1cm}\]
\item For $t>1$, $enum(t, i,\mbm{v},\mbm{w})$ is the formula 
\[
 \bigvee_{\min\{s_1,\ldots,s_{t-1}\}< s_{t}\, {\rm and } \, s_{t}= i \atop |\{s_1,\ldots,s_{t-1}\}|=t-1}\,\, \,( \bigwedge_{l=1}^{t}\bigwedge_{j=1}^{n} v_{s_l}\not= w_j\land \bigwedge_ {1\leq l<j\leq t}v_{s_l}\not=v_{s_j} \]\[\hfill \land \bigwedge_{ s < i \, {\rm and } \, s\not\in \{s_1,s_2,\ldots, s_{t-1}\}}\,\,\,( \bigvee_{j=1}^{n}
v_{s}= w_j\lor \bigvee_{1\leq j<t} v_{s} =v_{s_j})).\hspace*{1cm}\]
\end{description}

\begin{lemma}[A formula describing a permutation] Let $\S$ be a structure in ${\sf struc}_{\rm pred}^{({\rm m})}(I)$, $f$ be in ${\rm assgn}(\S)$, and $n,m\geq 1$. 
Let $P$ be a finite subset of $I$, let $q=|P|$, $\mbm{e} \in \{1,\ldots,q+1\}^n$, and ${\cal K} \in {\sf K}_{{\rm idx}_{\mbmss{e}}}$. Let $\{\beta_1,\ldots, \beta_{q+1}\}$ be the $P$-adequate partition of $I$ with $\beta_{q+1}=I\setminus P$, $\mbm{\mu}$ be any tuple in $\widetilde{\beta_{q+1}}^n$ with pairwise different components,
 and $\mbm{\xi}_{\mbmss{e}, {\cal K}}$ be the tuple $\mbm{\xi}_{\mbmss{\mu},\mbmss{e}, {\cal K}}^{I,P}$. Let
 $\mbm{\xi}$ be any tuple in $ I^n$ satisfying $\S\models_{f\langle{ B_1 \atop\beta_1} {\cdots\atop\cdots} { B_{q+1}\atop\beta_{q+1}}{\mbmty{x}\atop \mbmty{\xi}}\rangle}F_{\mbmss{e}, {\cal K}}(\mbm{x})$, $\pi_{\mbmss{\xi}}\in \G^I_{\sf 1}(P)$, and $\mbm{\eta}_0\in I^m$. Then,
\[\S\models_{f'} swap_{\mbmss{e}}(\mbm{x}_0,\mbm{y}_0,\mbm{x},\mbm{y})\] holds for $f'=f\langle {\mbmss{x}_0\atop\mbmss{\xi}_{\mbmty{e}, {\cal K}}}{\mbmss{y}_0\atop\mbmss{\eta}_0}{\mbmss{x}\atop\mbmss{\xi}}{\mbmss{y}\atop\pi_{\mbmty{\xi}}(\mbmss{\eta}_0)}\rangle $. 
\end{lemma}
The latter means that 
$f'=f\langle {\mbmss{x}_0\atop\mbmss{\xi}_{\mbmty{e}, {\cal K}}}{\mbmss{y}_0\atop\mbmss{\eta}_0}{\mbmss{x}\atop\mbmss{\xi}}{\mbmss{y}\atop\mbmss{\eta}}\rangle $ holds for $\mbm{\eta}=\pi_{\mbmss{\xi}} (\mbm{\eta}_0)$.

Now, the preparation for proving Proposition \ref{ChoiceInS1} in the next section is closed. 

\subsection{The Ackermann axioms in the basic Fraenkel model II}

 Let us recall that $\S_0$ is the Henkin-Asser structure $\S(\bbbn,\G_{\sf 1} ^\bbbn,\I_{\sf 0}^\bbbn)$.

\begin{proposition}\label{ChoiceInS1} For any $n,m\geq 1$, any $H$ in ${\cal L}_{\mbmss{x},D}^{(2)}$, and any $f\in {\rm assgn}(\S_0)$, there holds
\[\S_0\models choice_h^{n,m}(H).\]
\end{proposition}
{\bf Proof.} Let $\S$ be the basic structure $(J_n)_{n\geq 0}$ with the domains $J_0=\bbbn$ and $J_n=J_n(\bbbn,\G_{\sf 1} ^\bbbn,\I_{\sf 0}^\bbbn)$ for $n\geq 1$ and let $\G=\G_{\sf 1} ^\bbbn$. Let $f\in {\rm assgn}(\S)$, $n,m\geq 1$, and $H(\mbm{x},D)$ be a formula in ${\cal L}_{\mbmss{x},D}^{(2)}$ such that $\S\models_f \forall \mbm{x} \exists D H(\mbm{x},D) $ holds.
Let $P$ be a finite subset of $\bbbn$ such that \[\S_{f\langle{\mbmty{x}\atop \mbmty{\xi}}{D\atop\delta}\rangle}(H(\mbm{x},D))=\S_{f\langle{\mbmty{x}\atop \pi(\mbmty{\xi})}{D\atop\delta^\pi}\rangle}(H(\mbm{x},D))\] holds for all lists $\mbm{\xi}\in \bbbn^n$, all predicates $\delta\in J_{m}$, and all $\pi\in \G(P)$. Such a support $P$ exists by Lemma \ref{EndlichesPfuerH}. Let $\{ \beta_1,\ldots,\beta_{q+1} \}$ be the $P$-adequate partition of $I$ with $\beta_{q+1}=\bbbn \setminus P$. In the following, let $\mbm{\mu}$ be the tuple $(\mu_1, \ldots, \mu_n)$ defined by $\mu_1=\min (\bbbn \setminus P)$ and by $\mu_i=\min ((\bbbn \setminus P)\setminus \{\mu_1,\ldots,\mu_{i-1}\})$ for $i\in\{2,\ldots,n\}$ if $n\geq 2$ and let $P_{\mbmss{\mu}}=\{\mu_1,\ldots,\mu_n\}$.
 Let $B_1, \ldots, B_{q+1}$ be new variables for $1$-ary predicates that do not occur in $H(\mbm{x},D)$. For $\mbm{e}$, ${\cal K}$, and $F_{\mbmss{e}, {\cal K}}(\mbm{x})$ given as above (cf.\,\,Table \ref{overview_nota}) and for $\bar f=f\langle { B_1 \atop\beta_1} {\cdots\atop\cdots} { B_{q+1}\atop\beta_{q+1}}\rangle $, $\mbm{\xi}\in \bbbn^n$, $\delta \in J_m$, and each $\pi\in\G(P)$, we have 
\[\!\!\!\! \!\!\S_{\bar f\langle{\mbmty{x}\atop \mbmty{\xi}}{D\atop\delta}\rangle} (F_{\mbmss{e}, {\cal K}}(\mbm{x})\land H(\mbm{x},D))=\S_{\bar f\langle{\mbmty{x}\atop \pi( \mbmty{\xi})}{D\atop\delta^\pi}\rangle} (F_{\mbmss{e}, {\cal K}}(\mbm{x})\land H(\mbm{x},D))\] by Lemma \ref{Lemma_P_id} 
 (\ref{dieGs4Z}). 
Let $\alpha_0$ be the $n$-ary predicate $\alpha^{\bbbn,P}_{\mbmss{\mu}}$ defined as $\alpha^{I,P}_{\mbmss{\mu}}$ for any individual domain $I$ and considered in Lemma \ref{L_P_id} (see also Table \ref{overview_nota}). Moreover, for each $\mbm{e}\in \{1,\ldots,q+1\}^n$ and $ {\cal K} \in {\sf K}_{{\rm idx}_{\mbmss{e}}}$, let $\mbm{\xi}_{\mbmss{e}, {\cal K}}$ be the tuple $\mbm{\xi}_{\mbmss{\mu},\mbmss{e}, {\cal K}}^{\bbbn,P}$.
Thus, $\widetilde\alpha_0$ is the choice set $\{\mbm{\xi}_{\mbmss{e}, {\cal K}}\in \bbbn^n \mid \mbm{e} \in \{1,\ldots,q+1\}^n\land {\cal K} \in {\sf K}_{{\rm idx}_{\mbmss{e}}}\}$ considered in Lemma \ref{L_P_id}. $\widetilde\alpha_0$ is finite and, thus, $\alpha_0$ belongs to $J_n$. 
By assumption, $\S_{f\langle{\mbmty{x}\atop \mbmty{\xi}}\rangle} (\exists D H(\mbm{x},D))$ is $true$ for all $\mbm{\xi}\in \bbbn^n$. Hence, for any $\mbm{\xi}\in \alpha_0$ there is at least one predicate $\delta \in J_{m}$ such that 
\[\S\models_{f\langle{\mbmty{x}\atop \mbmty{\xi}} {D\atop\delta}\rangle}H(\mbm{x},D).\]
Consequently, we have a finite family $({\cal A}^H_{\mbmss{\xi}})_{\mbmss{\xi}\in \alpha_0} $ of non-empty sets $ {\cal A}^H_{\mbmss{\xi}}$ with $ {\cal A}^H_{\mbmss{\xi}}=\{\delta\in J_m\mid \S\models_{f\langle{\mbmty{x}\atop \mbmty{\xi}} {D\atop\delta}\rangle} H(\mbm{x},D)\}$ and \[\S\models_{f\langle
 {A\atop\alpha_0}\rangle} \forall \mbm{x}\exists D(A\mbm{x} \to H(\mbm{x},D)).\]
 This implies, by Lemma \ref{endlSigma}, the existence of a choice function $\varphi_s$ that can be defined analogously as in a proof of Lemma \ref{endlSigma} for $s=|\alpha_0|-1$ and the existence of an $(n+m)$-ary predicate (denoted by $\sigma_0$) that assigns $true$ to  the tuples in
 \[\{(\mbm{\xi}\,.\,\mbm{\eta})\in \bbbn^{n+m}\mid \mbm{\xi}\in \alpha_0 \,\,\&\,\, \mbm{\eta}\in\varphi_s ( {\cal A}^H_{\mbmss{\xi}})\}\] and $false$  to the other tuples. $\sigma_0$ can be defined  step by step as follows.
For each $ \mbm{\xi} \in \alpha_0$, let $\delta_{\mbmss{\xi}}$ be defined by $\delta_{\mbmss{\xi}}=\varphi_s( {\cal A}^H_{\mbmss{\xi}})$. 
Then, for all $\mbm{e} \in \{1,\ldots,q+1\}^n$ and ${\cal K} \in {\sf K}_{{\rm idx}_{\mbmss{e}}}$, $\delta_{\mbmss{\xi}_{\mbmty{e}, {\cal K}}} $ is --- by the definition of $ {\cal A}^H_{\mbmss{\xi}_{\mbmty{e}, {\cal K}}}$ --- in $\S$ and, therefore, there is some finite $P_{\delta_{\mbmty{\xi}_{\mbmty{e}, {\cal K}}}}\subseteq\bbbn$ such that ${\rm sym}_\G(\delta_{\mbmss{\xi}_{\mbmty{e}, {\cal K}}})\supseteq\G(P_{\delta_{\mbmty{\xi}_{\mbmty{e}, {\cal K}}}})$. Let \[\sigma_{(\mbmss{e}, {\cal K})}^{(0)}= \{(\mbm{\xi}_{\mbmss{e}, {\cal K}}\,.\,\mbm{\eta})\in \bbbn^{n+m}\mid\mbm{\eta}\in \delta_{\mbmss{\xi}_{\mbmty{e}, {\cal K}}}\}.\] 
Since ${\rm sym}_{\G}(\sigma_{(\mbmss{e}, {\cal K})}^{(0)})\supseteq \G(P\cup P_{\mbmss{\mu}}\cup P_{\delta_{\mbmty{\xi}_{\mbmty{e}, {\cal K}}}})$ follows from
 ${\rm sym}_{\G}(\{\mbm{\xi}_{\mbmss{e}, {\cal K}}\})\supseteq \G(P\cup P_{\mbmss{\mu}})$ (cf.\,\,Lemma \ref{L_P_id}), the predicate $\sigma_{(\mbmss{e}, {\cal K})}^{(0)}$ belongs to $J_{n+m}$.
Let the predicate $\sigma_0$ be given by 
\[\sigma_0= \bigcup_{\mbmss{e}\in \{1,\ldots,q+1\}^n}\,\bigcup_{{\cal K} \in {\sf K}_{ {\rm idx}_ {\mbmty{e}}}} \sigma_{(\mbmss{e}, {\cal K})}^{(0)}.\]
Since $\sigma_0$ is the predicate $\alpha _{\S,G _0 (\mbmss{x}, \mbmss{y}), (\mbmss{x}\,.\, \mbmss{y}), \tilde f}$ that can be defined by $ G_0 (\mbm{x}, \mbm{y})= \bigvee_{i=0,\ldots,s-1}R_i \mbm{x}\mbm{y} $ and $\tilde f= f\langle{R_0\atop \rho_0}{\cdots\atop\cdots}{R_{s-1}\atop \rho_{s-1}}\rangle$ if we assume that $ \{\rho_0,\ldots,\rho_{s-1}\}=\{\sigma_{(\mbmss{e}, {\cal K})}^{(0)}\mid \mbm{e}\in \{1,\ldots,q+1\}^n\,\,\&\,\,{\cal K} \in {\sf K}_{ {\rm idx}_ {\mbmty{e}}}\} $, it belongs to $J_{n+m}$. By Proposition 2.21, the relationship ${\rm sym}_{\G}(\sigma_0)\supseteq \bigcap_{\mbmss{e}\in \{1,\ldots,q+1\}^n}\,\bigcap_{{\cal K} \in {\sf K}_{ {\rm idx}_ {\mbmty{e}}}} \G(P\cup P_{\mbmss{\mu}}\cup P_{\delta_{\mbmty{\xi}_{\mbmty{e}, {\cal K}}}})$ is satisfied. Let $P_0$ be the finite set $ \bigcup_{\mbmss{e}\in \{1,\ldots,q+1\}^n}\,\bigcup_{{\cal K} \in {\sf K}_{ {\rm idx}_ {\mbmty{e}}}}P_{\delta_{\mbmty{\xi}_{\mbmty{e}, {\cal K}}}}$. Thus, we have ${\rm sym}_{\G}(\sigma_0)\supseteq\G(P\cup P_{\mbmss{\mu}}\cup P_0 )$ and the set $P\cup P_{\mbmss{\mu}}\cup P_0$ is a finite individual support for $\sigma_0$.
By definition of $\sigma_0$, there holds
 \[\S\models_{f\langle{A\atop\alpha_0}{S\atop \sigma_0} \rangle}\forall \mbm{x} \exists D ( A\mbm{x}\to \forall\mbm{y}(D \mbm{y} \leftrightarrow S\mbm{x}\mbm{y}) \land H(\mbm{x}, D )).\] 
Now, let $\mbm{e} \in \{1,\ldots,q+1\}^n$ and ${\cal K} \in {\sf K}_{{\rm idx}_{\mbmss{e}}}$. 
Let $\pi_{\mbmss{\xi}}$ be the permutation defined for all $\mbm{\xi}$ in Section \ref{DefinPerm} and let \[\sigma_{\mbmss{e}, {\cal K}}= 
\{(\mbm{\xi}\,.\,\mbm{\eta})\in \bbbn^{n+m}\mid \mbm{\eta}\in (\delta_{\mbmss{\xi}_{\mbmty{e}, {\cal K}}})^{\pi_{\mbmty{\xi}}}\}\] which means that the tuple $(\mbm{\xi}\,.\,\mbm{\eta})$ belongs to $\sigma_{\mbmss{e}, {\cal K}}$ iff there is an $\mbm{\eta}_0\in \delta_{\mbmss{\xi}_{\mbmty{e},{\cal K}}}$ such that $\mbm{\eta}=\pi_{\mbmss{\xi}}(\mbm{\eta}_0)$.
Let $\bar f=f\langle{ B_1 \atop\beta_1} {\cdots\atop\cdots} { B_{q+1}\atop\beta_{q+1}} \rangle $ and $f'=\bar f\langle {D\atop \delta_{\mbmty{\xi}_{\mbmty{e}, {\cal K}}}} {\mbmss{x}_0\atop\mbmss{\xi}_{\mbmty{e}, {\cal K}}} \rangle $. For 
\[G _{\mbmss{e}, {\cal K}} (\mbm{x}, \mbm{y} )= F_{\mbmss{e}, {\cal K}}(\mbm{x})\land 
\exists\mbm{y}_0(D\mbm{y}_0\land swap_{\mbmss{e}}(\mbm{x}_0, \mbm{y}_0,\mbm{x},\mbm{y})),\]
we get $\sigma_{\mbmss{e}, {\cal K}} =\alpha _{\S,G _{\mbmty{e}, {\cal K}} (\mbmss{x}, \mbmss{y} ), (\mbmss{x}\,.\, \mbmss{y}), f'}$.
 Consequently, $\sigma_{\mbmss{e}, {\cal K}} \in J_{n+m}$. By Proposition 2.21, ${\rm sym}_{\G}(\sigma_{\mbmss{e}, {\cal K}})\supseteq \G(P\cup P_{\mbmss{\mu}}\cup P_{\delta_{\mbmty{e}, {\cal K}}})$ holds. Let $f''= f'\langle{S_0\atop \sigma_0} \rangle$. Because of $\S\models_{f''} D\mbm{y}\leftrightarrow S_0 \mbm{x} _0\mbm{y}$, we have the equation $\S_{f'} ( G _{\mbmss{e}, {\cal K}} (\mbm{x}, \mbm{y} ))= \S_{f''}( G _{\mbmss{e}, {\cal K}} '(\mbm{x}, \mbm{y}) )$ for
\[G _{\mbmss{e}, {\cal K}} '(\mbm{x}, \mbm{y} )=F_{\mbmss{e}, {\cal K}}(\mbm{x}) \land\exists \mbm{x} _0 \exists \mbm{y}_0 ( F_{\mbmss{e}, {\cal K}}(\mbm{x}_0) \land S_0\mbm{x}_0\mbm{y}_0 \land swap_{\mbmss{e}}(\mbm{x}_{0},\mbm{y}_0,\mbm{x},\mbm{y})
 )\] 
which implies
\[\sigma_{\mbmss{e}, {\cal K}}=\alpha _{\S,G_{\mbmty{e}, {\cal K}} (\mbmss{x}, \mbmss{y} ), (\mbmss{x}\,.\, \mbmss{y}), f'} =\alpha _{\S,G' _{\mbmty{e}, {\cal K}} (\mbmss{x}, \mbmss{y} ), (\mbmss{x}\,.\, \mbmss{y}), f''}.\]
Let $\mbm{\xi}$ satisfy $\S\models_{\bar f\langle{\mbmty{x}\atop\mbmty{\xi}}\rangle} F_{\mbmss{e}, {\cal K}}(\mbm{x})$ and let $\delta_{\mbmss{\xi}}$ be the predicate $(\delta_{\mbmss{\xi}_{\mbmty{e}, {\cal K}}})^{\pi_{\mbmty{\xi}}} $. Then, there holds $\{\mbm{\eta}\mid (\mbm{\xi}\,.\,\mbm{\eta}) \in\sigma_{\mbmss{e}, {\cal K}}\}= \delta_{\mbmss{\xi}}$. Since ${\rm sym}_{\G}(\delta_{\mbmty{\xi}_{\mbmty{e}, {\cal K}}})\supseteq \G(P_{\delta_{\mbmty{\xi}_{\mbmty{e}, {\cal K}}}})$ implies ${\rm sym}_{\G}(\delta_{\mbmty{\xi}})\supseteq \G(\pi_{\mbmss{\xi}}(P_{\delta_{\mbmty{\xi}_{\mbmty{e}, {\cal K}}}}))$, $\delta_{\mbmss{\xi}}$ belongs to $\S$, too. 
By Lemma \ref{EndlichesPfuerH} and by the definition of $P$, for any $\pi\in \G(P)$, we have
\[\S_{f\langle {\mbmty{x}\atop \mbmty{\xi}_{\mbmty{e}, {\cal K}}}{D\atop\delta_{\mbmty{\xi}_{\mbmty{e}, {\cal K}}}}\rangle}(H(\mbm{x},D))=\S_{f\langle {\mbmty{x}\atop \pi(\mbmty{\xi}_{\mbmty{e}, {\cal K})}}{D\atop(\delta_{\mbmty{\xi}_{\mbmty{e}, {\cal K}}})^{\pi}}\rangle}(H(\mbm{x},D))=\S_{f\langle{\mbmty{x}\atop \mbmty{\xi}} 
 {D\atop\delta_{\mbmty{\xi}}}\rangle}(H(\mbm{x},D))\] which means $\delta_{\mbmss{\xi}} $ is in ${\cal A}_{\mbmss{\xi}}^H$ (defined in Lemma \ref {Zhg_CG_sigma}) and
\[\S\models_ {\bar f\langle {S\atop \sigma_{\mbmty{e}, {\cal K}}}\rangle}\forall \mbm{x} ( F_{\mbmty{e}, {\cal K}}(\mbm{x})\to \exists D (\forall\mbm{y}(D \mbm{y} \leftrightarrow S\mbm{x}\mbm{y}) \land H(\mbm{x}, D ))).\] 
On the other hand, the predicate
 \[\sigma= \bigcup_{\mbmss{e}\in \{1,\ldots,q+1\}^n}\,\,\bigcup_{{\cal K} \in {\sf K}_{{\rm idx}_{\mbmty{e}}}}\sigma_{\mbmss{e}, {\cal K}}\] that can be also described by 
$\sigma =\alpha _{\S,G (\mbmss{x}, \mbmss{y}), (\mbmss{x}\,.\, \mbmss{y}), f''}$
for 
\[G(\mbm{x}, \mbm{y})= \bigvee_{\mbmss{e}\in\{1,\ldots, q+1\}^n}\bigvee_{{\cal K}\in {\sf idx}_{\mbmty{e}}} G _{\mbmss{e}, {\cal K}} '(\mbm{x}, \mbm{y} )\]
 --- which implies ${\rm sym}_{\G}(\sigma)\supseteq \G(P\cup P{\mbm{\mu}}\cup P_0)$ --- belongs to $J_{n+m}$ and satisfies
\[\S\models_{f\langle{S\atop \sigma} \rangle}\forall \mbm{x} \exists D ( \forall\mbm{y}(D \mbm{y} \leftrightarrow S\mbm{x}\mbm{y}) \land H(\mbm{x}, D )). \] 
This means that the consequent of $choice_h^{n,m}(H)$ follows from the antecedent of $choice_h^{n,m}(H)$. \qed

\begin{corollary}\label{choiceInS1} For any $n, m\geq 1$ and any $H$ in ${\cal L}_{\mbmss{x},D}^{(2)}$, there holds 
\[\S_0\models choice^{n,m}(H).\]
\end{corollary}

A further conclusion can be derived from Corollary \ref{choiceInS1} and so on. 
\begin{corollary}\label{ACnmInS1} For any $n, m\geq 1$, there holds 
\[\S_0\models AC^{n,m}.\]
\end{corollary}

The following proposition results from Proposition \ref{ChoiceInS1}, Proposition \ref {PropertiesWO}, and Proposition \ref{WOnichtInS1}. 

\begin{proposition}\label{WO_unabh_choice} There is a Henkin-Asser structure that is a model of 
\[^{h}ax^{(2)} \cup choice_h^{n,m} \cup \{\neg WO^1\}\]
 for any $n\geq 1$ and $m\geq 1$.
\end{proposition}

\section{A summary}\label{Summary}

\setcounter{satz}{0}

 In \,{\rm HPL}, $AC_*^{1,1}$ is weaker than $AC^{n,m}$ for any $n,m\geq 1$. Proposition \ref{WO_unabh_choice} answers a question of Paolo Mancosu and Stewart Shapiro (2018). The main results of this section describe the following independence relations for HPL. 

By Proposition \ref{WO_unabh_choice} there is a Henkin-Asser structure that is a model of $\{\neg WO^{1}\}\cup choice_h$. Moreover, if we use ZFC as metatheory and ZFC is consistent, then every standard structure $\S_1^I$ is a model of $WO^{1}$ by \cite{Zerm04} and thus a model of $\{ WO^{1}\}\cup choice_h^{n,m}$. The latter follows from Lemma \ref{Zhg_CG_sigma} and so on.

\begin{theorem}\label{HPL_ind_3}$WO^1 $ is {\rm HPL}-independent of $ choice_h$.
\end{theorem}

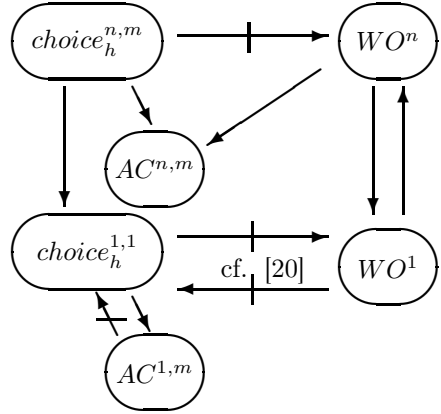
\begin{figure}[h] \unitlength1cm
\hspace{3cm}\begin{picture}(12,5.2) \thicklines
\put(0,-2.3){
\put(5.05,6.9){\oval(1.3,1)\makebox(0,0){$WO^{n}$}}
\put(5.25,4.6){\vector(0,1){1.7}}
\put(4.85,6.3){\vector(0,-1){1.7}}
\put(2.25,6.98){\vector(1,0){2}}
\put(3.2,6.75){\line(0,1){0.4}}
\put(4.15,6.52){\vector(-3,-2){1.5}}
\put(0,-3)
{
\put(5.05,6.9){\oval(1.3,1)\makebox(0,0){$WO^{1}$}}
\put(2.25,7.3){\vector(1,0){2}}
\put(4.25,6,6){\vector(-1,0){2}}
\put(3.25,6.4){\line(0,1){0.4}}
\put(3.25,7.1){\line(0,1){0.4}}
\put(3.3,6.85){\makebox(0,0){\, cf. \,\cite{Siskind}}}
}

\put(-3,2.0)
{
\put(4.95,3.2){\oval(1.3,1)\makebox(0,0){$AC^{n,m}$}}
\put(4.95,0.5){\oval(1.3,1)\makebox(0,0){$AC^{1,m}$}}
}
\put(1.75,-0.7)
{
\put(-0.7,7.6){\oval(2,1)\makebox(0,0){$choice_h^{n,m}$}}
\put(-1,7){\vector(0,-1){1.6}}
\put(-0.1,7){\vector(1,-2){0.25}}
\put(-0.7,4.8){\oval(2,1)\makebox(0,0){$choice_h^{1,1}$}}
\put(-0.1,4.2){\vector(1,-2){0.23}}
\put(-0.58,3.9){\line(1,0){0.4}}
\put(-0.3,3.7){\vector(-1,2){0.28}}
}
}
\end{picture}
\end{figure}

Note, that we also have the following theorems.
\begin{theorem}\label{HPL_ind_1} $choice_*^1$ is {\rm HPL}-independent of $AC^{1,1}$. 
\end{theorem}
\begin{theorem} \label{HPL_ind_2} $choice_h^{1,1}$ is {\rm HPL}-independent of $AC^{1,1}$. \end{theorem}

\end{document}